\documentclass[twocolumn,10pt]{article}

\usepackage{amsmath,amssymb,amsfonts}
\usepackage{mathptmx}
\usepackage[scaled=.81]{helvet}
\usepackage{courier}
\usepackage{newpxmath}
\usepackage{graphicx}
\iffalse
\usepackage{algorithmic}
\usepackage{textcomp}
\usepackage[dvipsnames,svgnames,x11names]{xcolor}
\usepackage{titling}
\pretitle{\begin{flushleft}\Large\sffamily\color{NavyBlue}}
\posttitle{\par\end{flushleft}\vskip-1em}
\preauthor{\begin{flushleft}\large\sffamily\color{darkgray}}
\postauthor{\end{flushleft}\vskip-1em}
\predate{\begin{flushleft}\normalsize\scshape\color{DarkGreen}}
\postdate{\par\end{flushleft}\vskip-1em}
\usepackage[bf,compact,small]{titlesec}
\titleformat*{\section}{\sffamily\bfseries\color{NavyBlue}}
\titleformat*{\subsection}{\sffamily\bfseries}
\titleformat*{\subsubsection}{\sffamily\bfseries}
\titleformat*{\paragraph}{\sffamily\bfseries}
\usepackage[margin=9ex]{geometry}
\else
\usepackage{savetrees}
\fi

\DeclareMathOperator{\nnzp}{nnzpattern}

\newcommand{\IEEEauthorblockN}[1]{#1}
\newcommand{\IEEEauthorblockA}[1]{\\#1}
\newenvironment{IEEEkeywords}{\textbf{Keywords:}}{}

\begin{document}

\newsavebox{\firstthnx}
\NewCommandCopy\origthnx\thanks
\renewcommand\thanks[1]{\sbox{\firstthnx}{#1}}

\title{Performance and Numerical Aspects of Decompositional
Factorizations with FP64 Floating-Point Emulation in INT8
\thanks{Research was sponsored by the Department of the Air Force Artificial Intelligence Accelerator and was accomplished under Cooperative Agreement Number FA8750-19-2-1000. The views and conclusions contained in this document are those of the authors and should not be interpreted as representing the official policies, either expressed or implied, of the Department of the Air Force or the U.S. Government. The U.S. Government is authorized to reproduce and distribute reprints for Government purposes notwithstanding any copyright notation herein.}
}


\author{\IEEEauthorblockN{Piotr Luszczek, Vijay Gadepally, LaToya
Anderson, William Arcand, David Bestor, William Bergeron, \\
Alex Bonn, Daniel J. Burrill, Chansup Byun, Michael Houle, Matthew
Hubbell, Hayden Jananthan, \\ Michael Jones, Peter Michaleas,
Guillermo Morales, Julia Mullen, Andrew Prout,
Albert Reuther, \\ Antonio Rosa, Charles Yee, and Jeremy Kepner}
\IEEEauthorblockA{\textit{LLSC, MIT Lincoln Laboratory, Lexington, MA, USA}}}

\date{September 26, 2025\origthnx{Research was sponsored by the
Department of the Air Force Artificial Intelligence Accelerator and was
accomplished under Cooperative Agreement Number FA8750-19-2-1000. The
views and conclusions contained in this document are those of the
authors and should not be interpreted as representing the official
policies, either expressed or implied, of the Department of the Air
Force or the U.S. Government. The U.S. Government is authorized to
reproduce and distribute reprints for Government purposes
notwithstanding any copyright notation herein.}~${}^{,}$\origthnx{Use of this work is controlled by the %
human-to-human license listed in Exhibit 3 of %
https://doi.org/10.48550/arXiv.2306.09267}}

\maketitle
\usebox{\firstthnx}
\begin{abstract}
Mixing precisions for performance has been an ongoing trend as the
modern hardware accelerators started including new, and mostly
lower-precision, data formats. The advantage of using them is a great
potential of performance gain and energy savings. The disadvantage are
the numerical issues not present in the standard-mandated
floating-point formats. Split integer emulation of FP64 takes this to
an extreme with the computation performed only by fixed-point tensor
core units. We present the new issues the emulation faces for
practical cases involving dense linear solver. We show extensive
numerical tests indicating the effect of extended numerical range of
matrix entries. We also scaled the input sizes to study the performance
and numerical profiles on the NVIDIA Hopper GPUs.
\end{abstract}


\section{Introduction}
\label{sec:intro}

Mixed-precision methods~\cite{doi:10.1177/10943420211003313} continue to
proliferate in computational and data sciences due to the many hardware
accelerators supporting unconventional number
formats~\cite{reuther2019llscsuvey}. This spurred the development of
emulation techniques that deliver the accuracy of standard
floating-point types while using compute units that lack support for
essential features from the IEEE 754 specification.  These approaches
straddle the trade-off between the reduction of accuracy and the gain in
performance and energy efficiency. This
originates in the use of tensor or matrix units for low-precision
computations as opposed to vector and/or scalar units for the
high-precision equivalents. The performance benefit is often much larger
than the raw bit-count mount suggest: the 16-bit computations with FP16
or BFloat16 are more than 4 times as efficient as the 64-bit operations
on FP64 data. The two terms: accuracy and precision will be clarified
below. Here, we only point out that we focus on the use of
these terms in numerical analysis rather than in data science and
machine learning, which often have the flexibility of incorporating
lesser formats into their algorithmic frameworks as indicated by many
generative and large language models.

The methods that especially stand out among the FP64 emulation efforts
are those that use only the integer units, INT8 in particular. They come
to prominence due to the inflection point brought about by the introduction
of NVIDIA's Blackwell design. For the first time in the high end
compute line of GPUs, the FP64 performance decreased while all other formats continue the drastic pace of multi-fold
improvements. The previous generation chip, NVIDIA Hopper, is rated at
over 50 Tflop/s in FP64 and over 1500 Top/s in INT8. The hardware's
sparsity feature raises it over 3000 Tera-op/s. Both of
these ratings take advantage of the tensor core units. By comparison,
Blackwell is rated at 40 Tflop/s in FP64 because the tensor core unit
instructions are relegated to the fused multiply-add units.
At the same time, the INT8 tensor cores are rated at nearly 5 Peta-op/s
(or 10 Peta-op/s with sparsity). To summarize, the advantage of INT8 over
FP64 grows from 30-fold in Hopper to well over 100-fold in Blackwell.
This creates an opportunity to exploit this imbalance and emulate FP64
using the INT8 tensor core units.

In this paper, we investigate the numerical and performance aspects of
floating-point emulation. We focus specifically on the dense methods for
linear systems where calculations using FP64 are always in a great
need of hardware acceleration.

\section{Emulation of FP64 and INT8 Splitting}

The basic theory of emulation is based on the error-free
transformations~\cite{ozaki2012errorfree} and predates the existence of
fast tensor core units in hardware. The implementations of this scheme
came later targeting low-precision floating-point
formats~\cite{mukunoki2020dgemmtc} and the use of integer units
followed~\cite{ootomo2022fp32}. Without going into details, the first
step is to split the product of dense
matrices $A,B \in \mathbb{R}^{n \times n}$ into individual components
that are mapped onto separate compute units:

\begin{equation}
AB = (A_1 + A_2) (B_1 + B_2) = A_1 B_1 + A_1 B_2 + A_2 B_1 + A_2 B_2
\label{eqn:absplit}
\end{equation}
where matrices $A_1, A_2, B_1, B_2$ are derived from $A$ and $B$,
respectively, by splitting the 53 bits of FP64 significand (or mantissa)
to fit the target compute units. For example, splitting onto 32-bit
units, we would use a multiplier of $2^{32}$ so that
$\|A_1\| \approx 2^{32} \|A_2\|$ meaning that $A_1$ stores the higher
bits of the mantissa and $A_2$ stores the lower ones. For 8-bit
splitting, we need to use more components to account for a shorter
elements of the hardware tensor cores. In general, the $k$-way splitting
can be represented as

\begin{equation}
AB = \left( \sum_{i=1}^k A_i \right) \left( \sum_{i=1}^k B_i \right) = \sum_{i , j = 1 \ldots k} A_i B_j.
\label{eqn:absplitn}
\end{equation}

As multiplication produces extra bits on output (up to twice as many
bits as either of the inputs) they must be rounded which happens
automatically in FP64 hardware. For splitting schemes such
as~(\ref{eqn:absplit}), the rounding is induced by truncating the
summation early and not performing the dropped matrix products such as
$A_2 B_2$. Or for the general formula (\ref{eqn:absplitn}), the rounding
would set a threshold $t$ and drop the components of the sum for which
$i+j>t$. The additional details involving numerical aspects of
rounding are out of scope here.

The important aspect of split-scheme emulation is its lack of exponent
because all bits in the split form come from the original significands
(mantissae) thus drastic difference in number ranges of matrix elements
causes representation issues. Consider a trivial 2-by-2 matrix multiply
with identity
\begin{equation}
  A B \equiv
  \left[ \begin{array}{ c c } 2^{+53} & 2^{+53} \\ 2^{-53} & 2^{-53} \end{array} \right]
  \left[ \begin{array}{ c c } 1 & 0 \\ 0 & 1 \end{array} \right]
\label{eqn:range2x2}
\end{equation}
with a wildly changing numerical range between the 1\textsuperscript{st}
and 2\textsuperscript{nd} row of
matrix $A$.  Any fixed-width storage scheme such as the one used for the
split-scheme emulation of FP64 would round down the second row of $A$ to
$0$ due to its limited number of allocated space for the bits of split
matrices. In this exaggerated example the loss of accuracy due to the
dropped values seems small compared to the large norm of $A$ but in
practice this effect compounds for larger matrices as we show below
in our numerical experiments. Here, we only indicate how this problem
with numeric range could result from poor scaling or even a permutation of
rows and/or columns that cluster together matrix entries in a way that
causes truncation of mantissa bits when switching to the split integer
representation. Consider the original matrix product with scaling
applied to both arguments on both sides:

\begin{equation}
(L_A A R_A) (L_B B R_B)
\label{eqn:scaling}
\end{equation}
where we select the left and right scaling matrices to increase the
dynamic range of the entries in either $A$ or $B$. This could be
performed with trivial diagonal matrices that have a diagonal entries
varying widely in dynamic range just as the rows of matrix $A$ in
(\ref{eqn:range2x2}). Another possibility is to choose the left/right
matrices, say $L_A$ and $R_A$, to be a permutation matrices that group
together rows and columns of $A$ by increasing their dynamic range.
Either diagonal scaling or permutation matrices are cheap to apply and
invert so they are rarely considered to be an issue for calculations
in FP64.

The issues with split integer emulation presented so far were applied to
a generic matrix product. While they are concerning in general, they may
not be of practical importance for the common dense operations involving
matrix product. We next turn our attention to the case of dense linear
solvers where the numeric range issues could be masked for some
matrices, yet the problems persists if they are introduced through a
particular selection of matrix elements.

\section{Related Work}
\label{sec:related}

The foundational work that defined the error-free computations is
now commonly referred as the Ozaki Scheme~\cite{ozaki2012errorfree}.
In order to make it more widely applicable, the constant values
used for splitting were subsequently
generalized~\cite{ozaki2013generrorfree}. An implementation was later
revised to further improve the accuracy profile for the 2-slice
scenario~\cite{ozaki2015imperrorfree}.
To determine the best splitting in practice, a method of was proposed to
both test and validate the various
choices~\cite{ozaki2016errorfreevalid}. Unlike the original work, there
is also a potential to apply integer splitting to the case of sparse
matrix multiplication~\cite{ichimura2018spsmatvec}. The double precision
matrix-matrix multiply or DGEMM() emulation was the original application
but the use of FP16 Tensor Cores was also possible leveraging another
piece of computational hardware units inside modern
accelerators~\cite{mukunoki2020dgemmtc}. Another extension of the
integer splitting approach was to extend the precision of the emulated
format focusing on quadruple precision or FP128 for the matrix
multiplication operation~\cite{mukunoki2021matmul128}.
An alternative approach to error-free transformation used both the
Cauchy–Schwarz inequality and Karatsuba algorithm to arrive at a matrix
multiplication algorithm~\cite{lange2022altozaki}.
The derivation of the best splitting was further enhanced with an
iterative procedure~\cite{ozaki2022errapps}. An arbitrary length
precision was applied to both inner products and sparse matrix-vector
multiplications~\cite{mukunoki2022infprecdot}.
Improved accuracy enhancing efficiency for symmetric rank-k operations
were done in the context of eigenvalue and singular value
solvers~\cite{uchino2024eigsvd}. The recent efforts implemented DGEMM()
using emulation with INT8 Tensor Cores~\cite{ootomo2024dgemmint}.
Further acceleration of DGEMM() emulation was achieved by the means of
combing INT8 Tensor Cores with application of error-free computation
method~\cite{uchino2025ozakiperf}. Finally, for the case of unequal
sizes of integer splits, an accuracy extension was proposed to handle
the original error-free transformation~\cite{ozaki2025errfreeext}.

\section{Dense Numerical Solvers Using FP64 Emulation}

Decompositional factorizations are commonly used in computational and
data sciences as both efficient and numerically stable methods for
obtaining a solution of a system of linear equations. In its generic
form written as

\begin{equation}
Ax=b
\label{eqn:axb}
\end{equation}
where $A \in \mathbb{R}^{n \times n}$ and $x,b \in \mathbb{R}^n$. We are
concerned here with general non-singular matrices but it is possible our
methods would equally apply to symmetric cases either positive definite
or indefinite.
Given a two-by-two partitioning $A$ from (\ref{eqn:axb})

\begin{equation}
\left[ \begin{array}{ c c } A_{11} & A_{12} \\ A_{21} & A_{22} \end{array} \right]
\left[ \begin{array}{ c } x_1 \\ x_2 \end{array} \right]
=
\left[ \begin{array}{ c } b_1 \\ b_2 \end{array} \right]
\label{eqn:axb2x2}
\end{equation}
we can use the decompositional approach by factoring $A$ into a product
of triangular matrices producing the LU factorization as is the case for the
HPL benchmark~\cite{Dongarra2003}. In this process, all 3\textsuperscript{rd} order
flops are concentrated in the Schur complement:
\begin{equation}
A'_{22} \leftarrow A_{22} - A_{21} A^{-1}_{11} A_{12}
=
A_{22} - A_{21} (L_{11} U_{11})^{-1} A_{12}
\label{eqn:schur}
\end{equation}
which is a perfect target for fast matrix-matrix multiplication such as
the one implemented by the split-scheme emulation of FP64. The question
is then if there is a potential of adverse scaling from
(\ref{eqn:scaling}) to enter the factorization process thus forcing
truncation and increasing the errors beyond standard floating-point
roundoff. The basic observation is that the update formula
(\ref{eqn:schur}) has self-stabilizing property stemming from the
pivoting process that chooses the diagonal elements for scaling $A_{11}$
to remain numerically stable under mild assumptions about the
initial conditioning of $A$. Using the 2-by-2 partitioning of
(\ref{eqn:axb2x2}), partial pivoting (the most common pivoting scheme for
LU factorization) selects pivots from $[ A_{11} \: A_{21}]^T$ and produces
triangular $L_{11}$ that has $1$'s on the diagonal and entries of
magnitude less than $1$ below the diagonal. This leaves $U_{11}$ to be
the potential source of issues with large disparity in numeric range as
this upper triangular factor remains unscaled by the pivot value used
for $L_{11}$. However, the application of $U^{-1}_{11}$ to $A_{21}$
occurs inside the panel factorization (the LU factorization of a
rectangular portion of $A$) and it is done by a custom code that is
unlikely to benefit from increased performance of the split-scheme
emulation. As a result, it remains challenging to observe problematic
scaling inside the Schur complement, which leaves us to consider a
different type of matrices that may potentially create issues for the
split-scheme emulation.

In order to measure the numerical errors in the computed solution, we
will be using the \emph{scaled residual} defined as

\begin{equation}
\text{Scaled residual} =
\frac{ \|Ax-b\|_{\infty} }{ (\|A\|_{\infty} \|x\|_{\infty} + \|b\|_{\infty}) n \epsilon }
\label{eqn:sres}
\end{equation}
which needs to be below $16$ for a valid HPL run.~\cite{Dongarra2003}

\section{Parametric Family of Matrices with Challenging Numerical Range}
\label{sec:parawilk}

We turn now to a classic example of matrix showing exponential element
growth when using partial pivoting. It is commonly attributed to
Wilkinson~\cite{wilkinsonRoundingErrorsAlgebraic1963,wilkinsonAlgebraicEigenvalueProblem1965}
and it is simply an almost lower-triangular matrix with $1$'s on the
diagonal, $-1$'s below the diagonal, and an extra $1$'s above the
diagonal in the last column. A 5-by-5 example looks as following:

\begin{equation}
\text{Wilkinson}_5 =
\left[ \begin{array}{ r r r r r }
 1 &  0 &  0 &  0 & 1 \\
-1 & -1 &  0 &  0 & 1 \\
-1 & -1 &  1 &  0 & 1 \\
-1 & -1 & -1 &  1 & 1 \\
-1 & -1 & -1 & -1 & 1 \\
\end{array} \right].
\label{eqn:wilk}
\end{equation}

Factorizing the Wilkinson matrix of size $n$ using partial pivoting
produces $2^n$ element growth and thus it requires complete pivoting to
limit the pivot growth to a 3\textsuperscript{rd} degree
polynomial~\cite{trefethenAveragecaseStabilityGaussian1990}. This is
impractical because of how complete pivoting inhibits parallelism. But
the element growth of Wilkinson matrices is the property allowing us to
arbitrarily increase the numerical range of the matrix elements, which
is challenging to handle by the split-scheme emulation of FP64.

We first limit the element growth in the Wilkinson matrices by observing
that the inverse of the Turing matrix (no extra $1$'s in the last
column) has the largest element $2^{n-2}$~\cite{turing1948rounderrmtx}
in the bottom row and it is the consequence of the summation of all
previous column elements of the inverse:

\begin{equation}
\left[ \begin{array}{ r r r r r }
 1 &  0 &  0 &  0 & 0 \\
-1 &  1 &  0 &  0 & 0 \\
-1 & -1 &  1 &  0 & 0 \\
-1 & -1 & -1 &  1 & 0 \\
-1 & -1 & -1 & -1 & 1 \\
\end{array} \right]^{-1}=
\left[ \begin{array}{ r r r r r }
 1 &  0 &  0 &  0 & 0 \\
 1 &  1 &  0 &  0 & 0 \\
 2 &  1 &  1 &  0 & 0 \\
 4 &  2 &  1 &  1 & 0 \\
 8 &  4 &  2 &  1 & 1 \\
\end{array} \right].
\label{eqn:wilk0}
\end{equation}

This exponential growth can be mitigated by limiting how many $-1$'s are
below the diagonal and indicate this with parameter $d$. For $n=5$ and
$d=2$ we have

\begin{equation}
\left[ \begin{array}{ r r r r r }
 1 &  0 &  0 &  0 & 0 \\
-1 &  1 &  0 &  0 & 0 \\
-1 & -1 &  1 &  0 & 0 \\
 0 & -1 & -1 &  1 & 0 \\
 0 &  0 & -1 & -1 & 1 \\
\end{array} \right]^{-1}=
\left[ \begin{array}{ r r r r r }
 1 &  0 &  0 &  0 & 0 \\
 1 &  1 &  0 &  0 & 0 \\
 2 &  1 &  1 &  0 & 0 \\
 3 &  2 &  1 &  1 & 0 \\
 5 &  3 &  2 &  1 & 1 \\
\end{array} \right],
\label{eqn:wilkd}
\end{equation}
which leads us to observe that the numbers in columns form the
Fibonacci sequence: $F_n = F_{n-1} + F_{n-2}$. This is still an
exponentially growing sequence but it is also exponentially slower than
the original. In fact, we can regulate the base of the exponent by
adjusting the $d$ parameter because the growth of values in the inverse
matrix is determined by a generalized Fibonacci sequence $G^{(d)}_n$:

\begin{equation}
G^{(d)}_n = \sum_{i=1}^d G_{n-i}.
\label{eqn:genfib}
\end{equation}
It is easy to establish that $F_n \equiv G^{(2)}_n$ and $G^{(n)}_n=2^n$
when we set $G^{(n)}_0=1$ and $G^{(n)}_i = 2^{i-1}$ for $1 \le i \le n$.

The next characteristic of the Wilkinson growth is a long running
sequence of $1$'s that allows expansion of a long sequence with
geometric doubling. By using the $d$ parameter we are able to limit the
exponent base by choosing the right $G^{(d)}_n$ sequence. To limit the
sequence's length, consider an example with $n=5, d=5$ and a blocking
factor $b=2$:

\begin{equation}
\left[ \begin{array}{ r r r r r }
 1 &  0 &  1 &  0 & 1 \\
-1 &  1 &  1 &  0 & 1 \\
-1 & -1 &  1 &  0 & 1 \\
-1 & -1 & -1 &  1 & 1 \\
-1 & -1 & -1 & -1 & 1 \\
\end{array} \right]^{-1}=
\left[ \begin{array}{ r r r r r }
\frac{1}{2}&-\frac{1}{5}&-\frac{1}{5}&  0 & 0 \\
0 & \frac{1}{2}&-\frac{1}{2}&  0 & 0 \\
0 &  0 & \frac{1}{2}&-\frac{1}{5}&-\frac{1}{5}\\
0 &  0 &  0 & \frac{1}{2}&-\frac{1}{2}\\
\frac{1}{2}& \frac{1}{5}&-\frac{1}{5}& \frac{1}{5}&\frac{1}{5}\\
\end{array} \right],
\label{eqn:wilkb}
\end{equation}
where the blocking factor $b$ becomes the length of the sequence with
exponential growth in the matrix inverse.

Finally, we observe that the geometric sequences formed in our
parametrized form of the Wilkinson matrix depend on the initial value of
the in the upper part of the matrix. We thus replace those upper-part
$1$'s with the third and final parameter $\alpha$. We arrive at the
3-parameter formulation that we refer to as ParaWilk$_n(d,b,\alpha)$.

\section{Randomization for Parametric Matrices}
\label{sec:rand}

The common problem with parametrized matrix families is that they may
admit closed-form inverses and thus are problematic for data generators
in the context of scalable benchmarks.~\cite{fasi2021nopivlugen}
This leads us down the path of randomization, which is already used by
HPL~\cite{Dongarra2003} with the matrix elements drawn from a uniform
distribution ranging between $-1/2$ and $1/2$ denoted as $U(-1/2,1/2)$.
Such random distributions tend to generate matrices with favorable
numerical properties for solvers based on the standard floating-point
solvers.~\cite{luszczek2020scaldatagen} However, having matrix entries
with alternating signs creates a greater likelihood for cancellation of
growth in number range when performing the Schur complement
(\ref{eqn:schur}). To counteract this, we switch our matrix generator to
produce only positive random entries skewed towards larger values. In
particular we used $2U^2_n(0, 1)$ to produce the results in the numerical
experiments section. However, the randomization is not applied uniformly
in order to preserve the numerical pattern created by the ParaWilk
formula from Section \S\ref{sec:parawilk}. We use nonzero pattern
function
\begin{equation}
\nnzp(A) = \begin{cases}
1 & \text{if $a_{ij} \ne 0$}, \\
0 & \text{if $a_{ij} = 0$}.
\end{cases}
\label{eqn:nnzpattern}
\end{equation}
and only add the random entries where the ParaWilk is $0$:
\begin{equation}
  A \sim ( 1 - \text{ParaWilk}(d,b,\alpha) ) \odot 2U^2(0,1) + \text{ParaWilk}(d,b,\alpha)
\label{eqn:rndgen}
\end{equation}
where $\odot$ is Hadamard product or element-wise multiply.

\section{Numerical Properties of Linear Solver with INT8 Emulation}

\begin{table}[tb]
\caption{Scaled residual error for a solve of a system of linear
equations with a new data generation scheme
ParaWilk$_{256}(d=4,b=15,\alpha=1/2)$. Passing value is assumed to be
$16.0$ and requires at least 7 INT8 splits. The bottom row shows result
with hardware-native FP64.}
\centering
\begin{tabular}{ c r }
\hline
INT8 splits & Scaled residual \\
\hline
3 & 26464646.4755162261 \\
4 & 341348.2056036110 \\
5 & 3798.4563646804 \\
6 & 39.1393959359 \\
7 & 0.3245150561 \\
8 & 0.0030077355 \\
9 & 0.0001860455 \\
\hline
FP64 (cuBLAS) & 0.0002377248 \\
\hline
\end{tabular}
\label{tab:sresozk}
\end{table}

We begin the numerical results by showing in Tab.~\ref{tab:sresozk} the
increasing number of INT8 splits and the corresponding values of the
scaled residual given by (\ref{eqn:sres}): the main metric that
determines correctness of an HPL run~\cite{Dongarra2003}. For the
results in the table, the data generator in the tests produced the
following matrix: ParaWilk$_{256}(d=4, b=15, \alpha=1/2)$. There are a
few conclusions that can be reached from this basic test. Firstly, it is
possible to use INT8 splitting to emulate FP64 for this matrix. However,
it required as many as 7 splits to reach the passing value of the scaled
residual, which is prohibitive from the performance standpoint as shown
in the next section. Additionally, as many as 9 splits were needed to
reach the same value of the scaled residual as was achieved with
hardware-native FP64 which is shown at the bottom row of the table.
Thus, our parametric data generator can be used as a detector for using
FP64 emulation with integer splitting.

\begin{figure}[htb!]
\centering
\includegraphics[width=\linewidth]{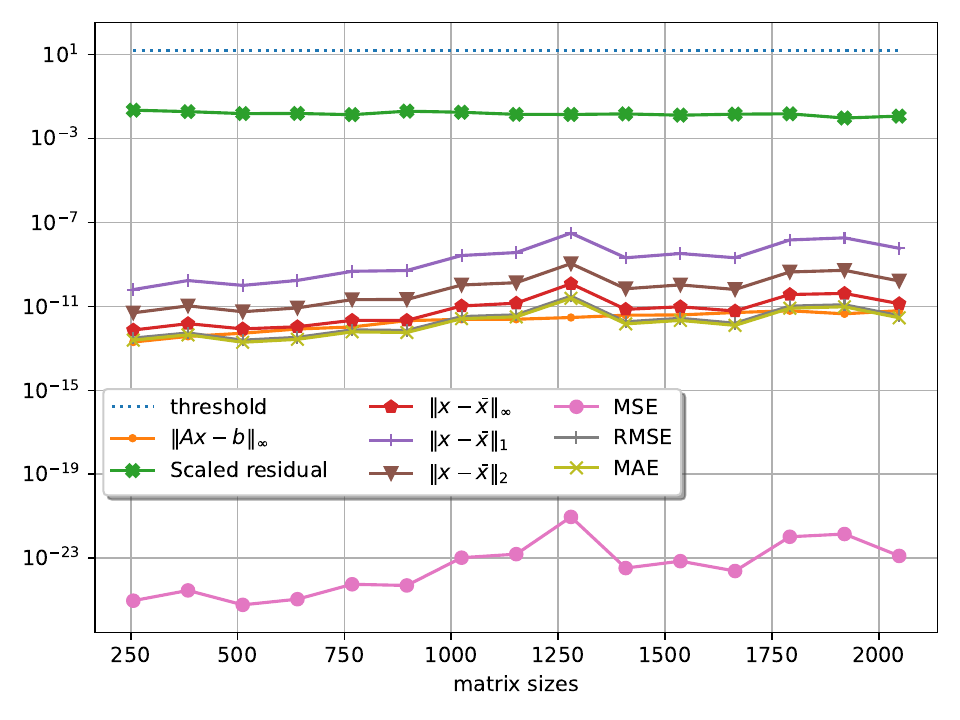}
\includegraphics[width=\linewidth]{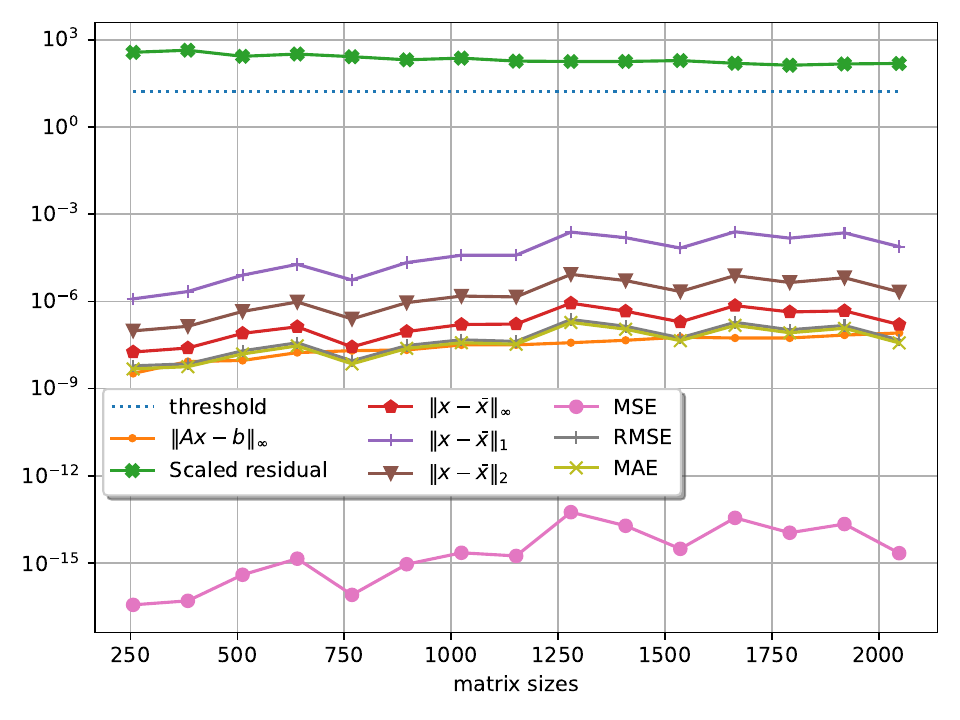}
\includegraphics[width=\linewidth]{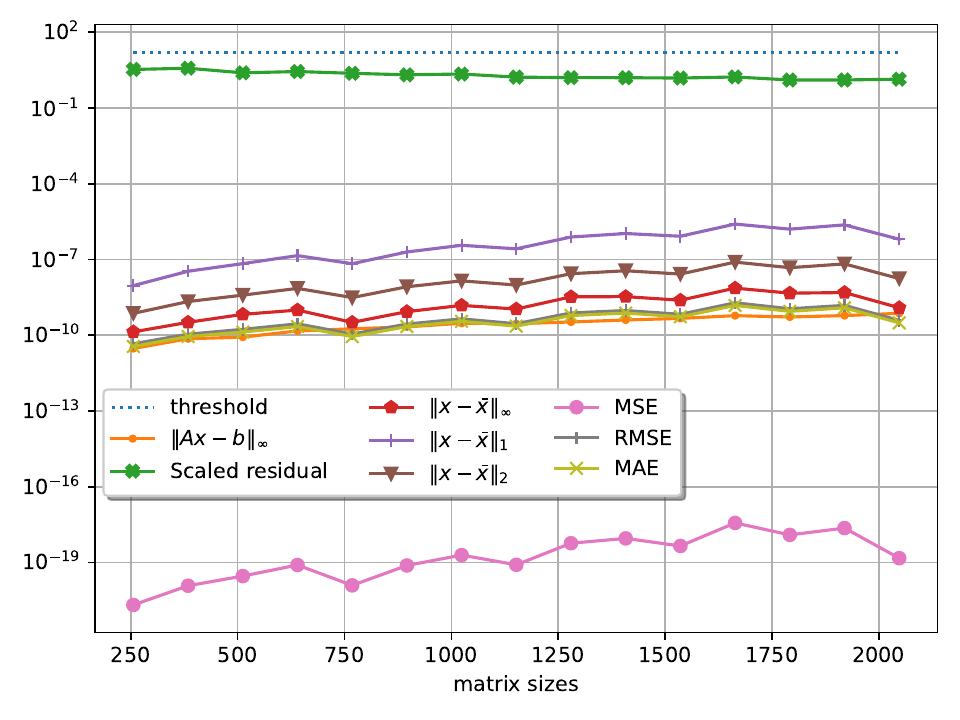}
\caption{Comparison of numerical and statistical metrics for the
solution of a system of linear equations of increasing sizes with
$A \sim U(-1/2,1/2)$. The results from using FP64 hardware in cuBLAS is
shown in the top. The bottom two plots show 6 and 7 INT8 splits,
respectively, and only the latter achieves the passing values of the
scaled residual: 16 (dotted line on all plots).}
\label{fig:errmetrics}
\end{figure}

Next we turn to more comprehensive numerical results shown in
Fig.~\ref{fig:errmetrics}, which intends to shed more light on the
scaling properties of our proposed parametric family of matrices as
their sizes increase. The figure shows a variety of numerical and
statistical metrics to indicate their behavior in the tested range of
matrix dimensions as the implementation of the LU factorization changes
from using the hardware-native FP64 at the top to using either 6 or
7 INT8 splits shown on two bottom plots, respectively. The matrix used in
these tests was generated with the standard HPL generator that used
uniform RNG centered around $0$: $A \sim U(-1/2, 1/2)$.

\begin{table}[tb]
\caption{Selected parameter values of the ParaWilk$_{n}(d, b, \alpha)$
family of matrices and the scaled residual value after LU factorization
and solve using the integer splitting emulation of FP64. We used
$\alpha=1$ for experiments in the table.}
\centering
\begin{tabular}{ c c c c c c c }
\hline
$n$ & \multicolumn{3}{c}{$6$ splits} & \multicolumn{3}{c}{$7$ splits} \\
\hline
    &     &     & Scaled   &     &     & Scaled \\
    & $d$ & $b$ & residual & $d$ & $b$ & residual \\
\hline
256 & 4 & 15 & 39.1 & 10 & 15 & 17.3 \\
384 & 4 & 18 & 19.3 & 11 & 12 & 18.9 \\
512 & 6 & 16 & 19.6 & 12 & 16 & 17.6 \\
640 & 7 & 22 & 21.7 & 13 & 16 & 24.1 \\
768 & 6 & 19 & 17.6 & 13 & 25 & 19.3 \\
896 & 8 & 15 & 40.3 & 14 & 28 & 19.0 \\
1024 & 8 & 23 & 22.4 & 15 & 19 & 16.3 \\
1152 & 8 & 15 & 40.5 & 16 & 18 & 18.2 \\
1280 & 8 & 16 & 25.8 & 16 & 30 & 18.9 \\
1408 & 9 & 16 & 22.2 & 18 & 22 & 24.6 \\
1536 & 9 & 28 & 27.5 & 17 & 23 & 26.3 \\
1664 & 11 & 19 & 40.6 & 17 & 26 & 16.7 \\
1792 & 10 & 23 & 16.5 & 19 & 23 & 61.6 \\
1920 & 10 & 15 & 17.0 & 18 & 20 & 38.8 \\
2048 & 11 & 23 & 20.7 & 19 & 26 & 19.3 \\
\hline
\end{tabular}
\label{tab:paramssres}
\end{table}

Finally, in Tab.~\ref{tab:paramssres} we show results from a scaling
experiments that indicate it is possible to further extend the numerical
accuracy effects in the LU factorization when choosing the right
parameters of our family of matrices given by the ParaWilk formulation
presented in Section \S\ref{sec:parawilk}. In order to limit the number
of experiments, we only present the first set of parameters we
encountered that caused the scaled residual exceed its maximum value of
$16$. Generally speaking, the values of $d$, or the number of elements
in the $G^{(d)}_n$ sequence (\ref{eqn:genfib}), are slowly monotonically
increasing with $n$ in a faster-than-logarithmic rate. The blocking
parameter $b$ is non-monotonic but its growth has the same order. This
is a positive indicator that the data generator is suitable to use with
larger matrix sizes.

\section{Performance Results}

\begin{figure}[htb]
\centering
\includegraphics[width=.95\linewidth]{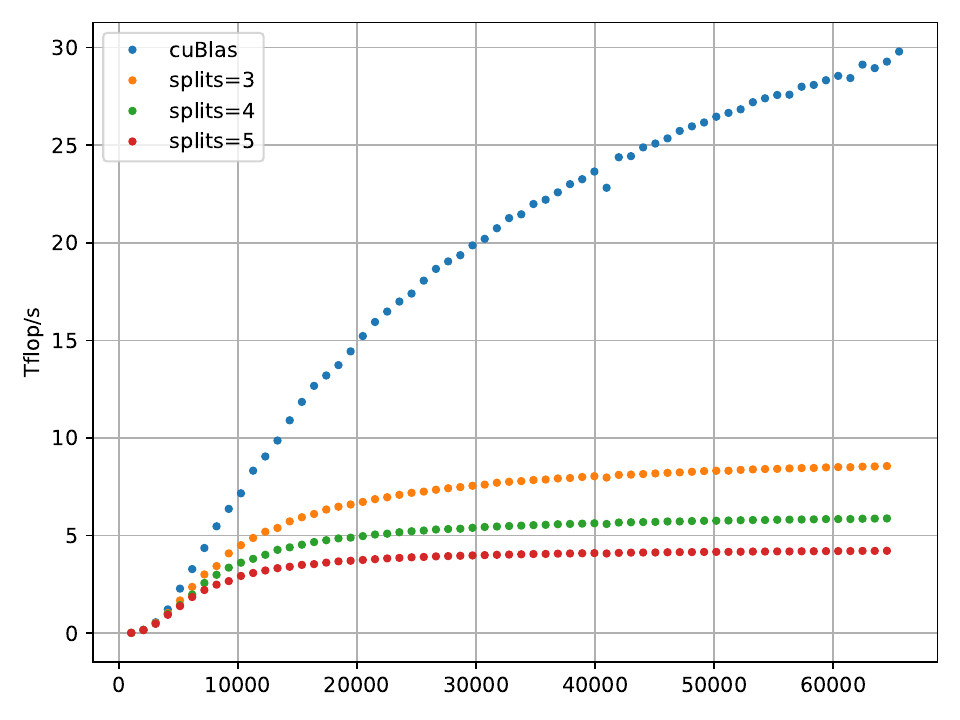}
\includegraphics[width=.95\linewidth]{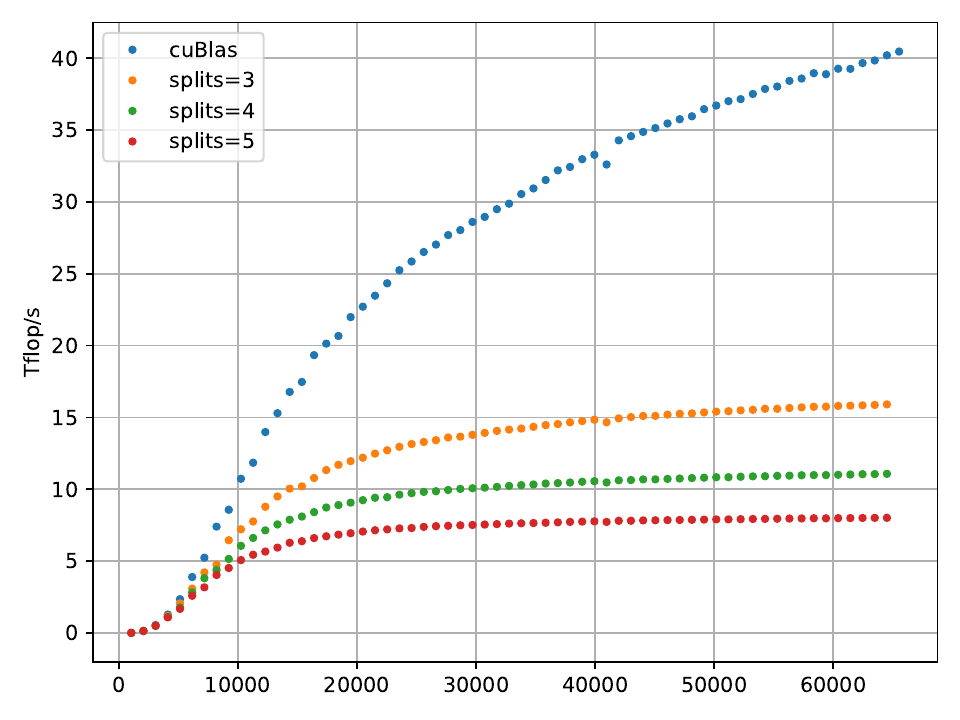}
\includegraphics[width=.95\linewidth]{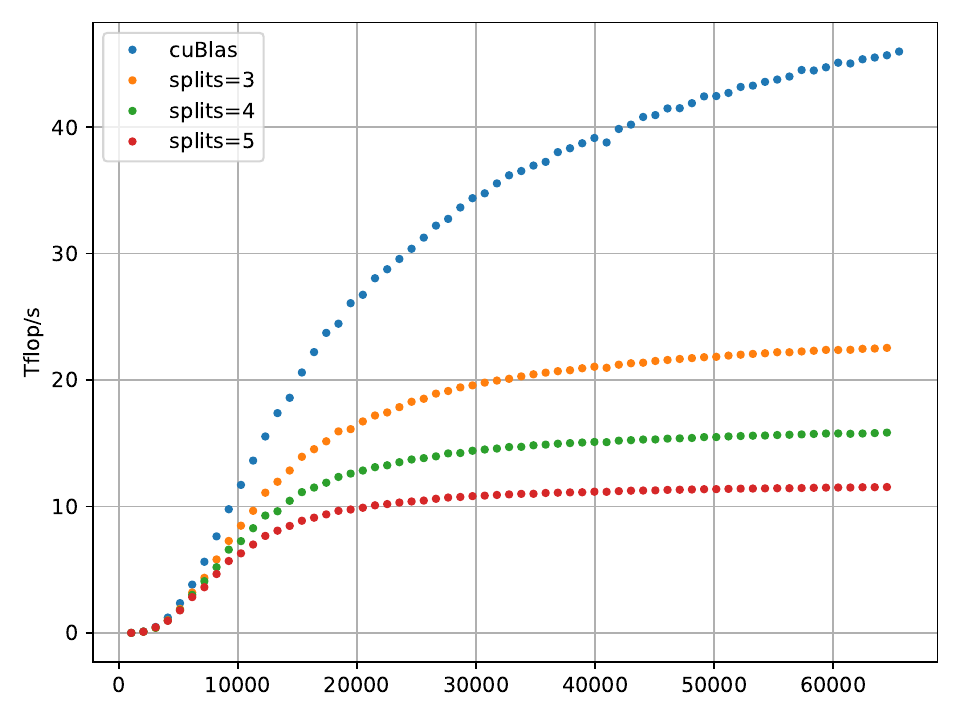}
\caption{Performance of cuBLAS in FP64 compared to split integer
implementation with varying number of splits and blocking factors:
256 (top), 512 (center), and 768 (bottom).}
\label{fig:ozkperf1}
\end{figure}

\begin{figure}[htb]
\centering
\includegraphics[width=.95\linewidth]{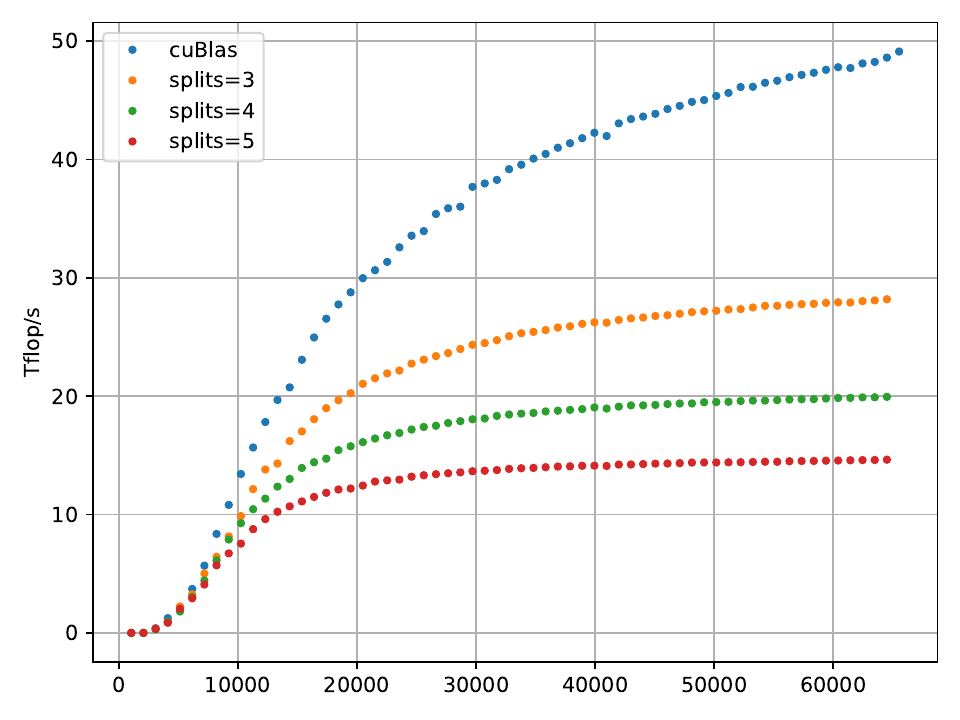}
\includegraphics[width=.95\linewidth]{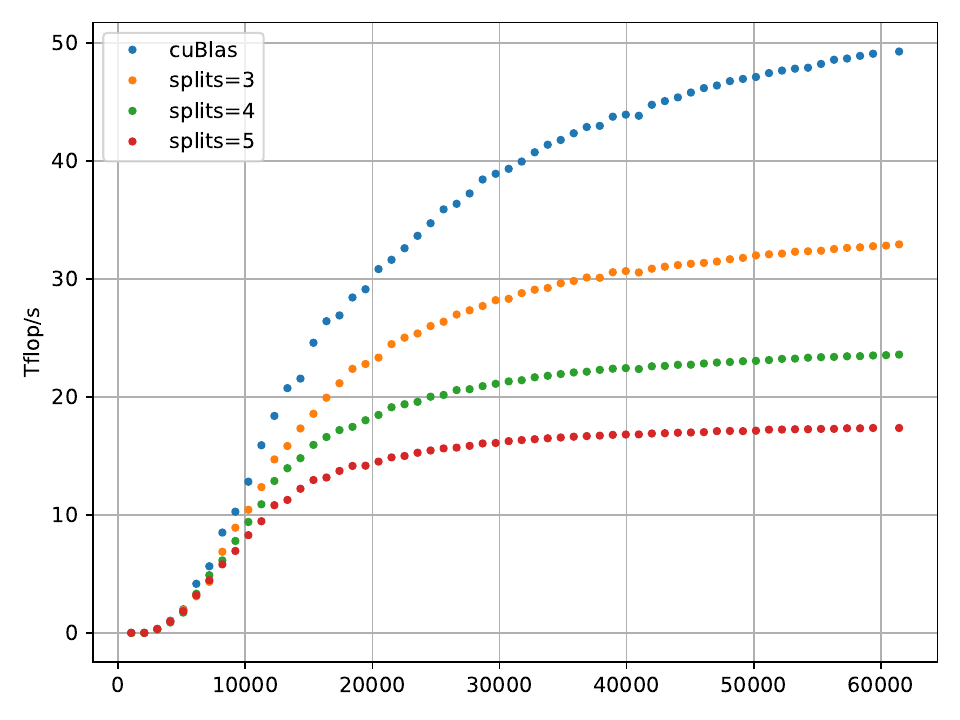}
\includegraphics[width=.95\linewidth]{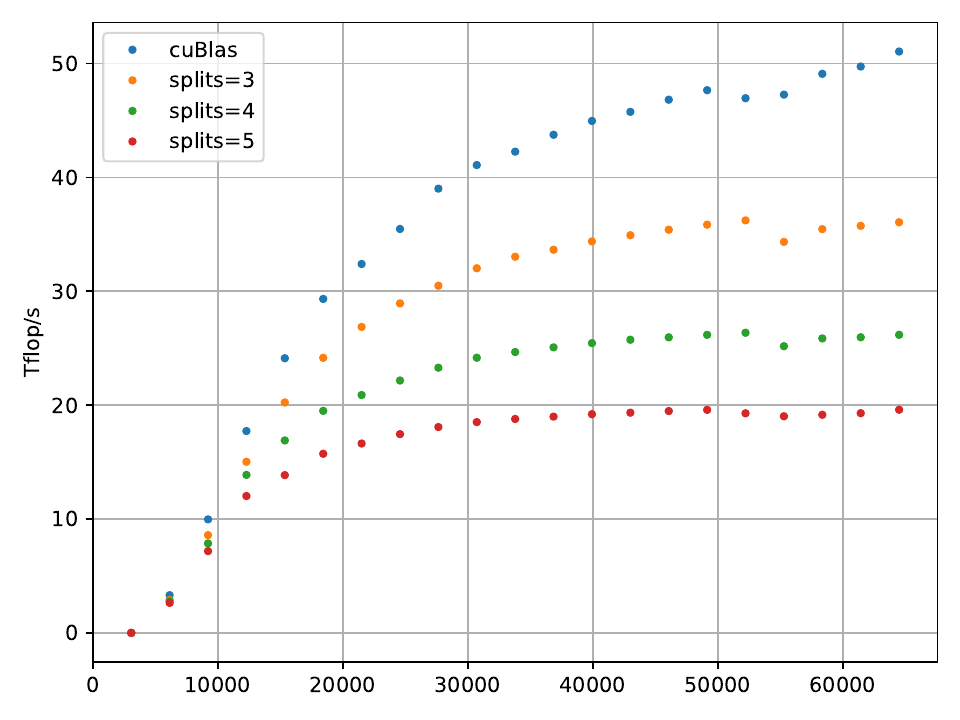}
\caption{Performance of cuBLAS in FP64 compared to split integer
implementation with varying number of splits and blocking factors:
1024 (top), 1280 (center), and 1536 (bottom).}
\label{fig:ozkperf2}
\end{figure}

\begin{figure}[htb]
\centering
\includegraphics[width=.95\linewidth]{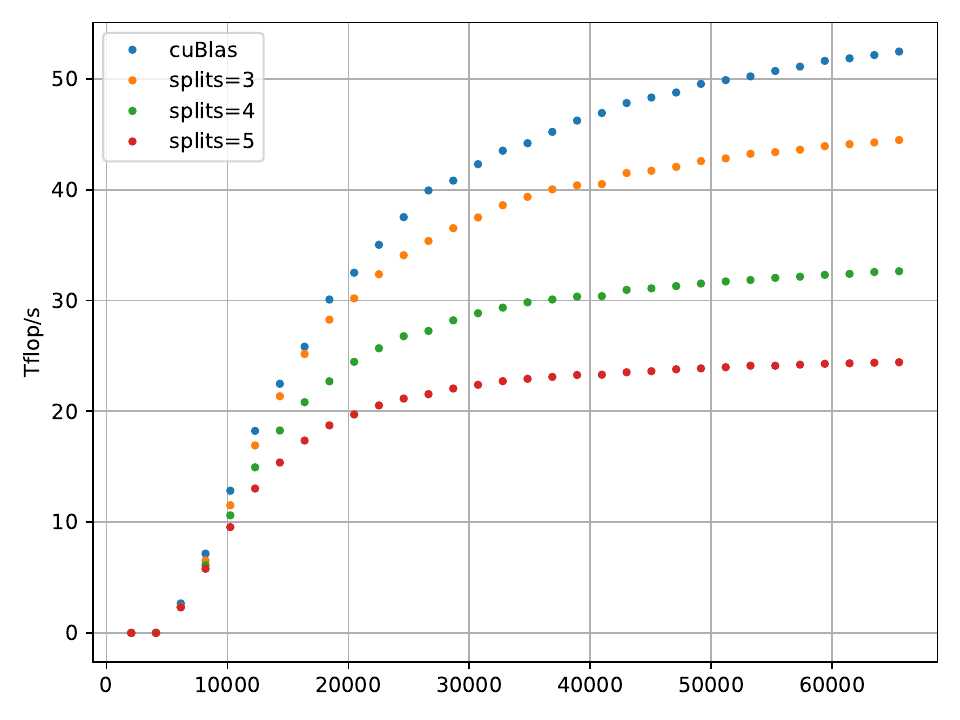}
\caption{Performance of cuBLAS in FP64 compared to split integer
implementation with varying number of splits and blocking factor of 2048.}
\label{fig:ozkperf3}
\end{figure}

In this section, we focus on the performance results to indicate the
benefits of INT8 split emulation of FP64 on a widely available NVIDIA
Hopper card.  In particular, we used PCIe version 4 model of this GPU
with 114 Streaming Multiprocessors and over 55 Tflop/s of 
FP64 performance at the Lincoln Laboratory Supercomputing
Center.~\cite{reuther2018interactive}

For our software implementation, we used the integer split emulation
code available in the Microsoft Github repository under handle
\textsf{enp1s0/ozIMMU} which is the basis for the official RIKEN
repository Microsoft Github located at
\textsf{RIKEN-RCCS/accelerator\_for\_ozIMMU}.

In order to focus the obtained performance values on the performance
obtained from the subsequent Schur complements (\ref{eqn:schur}), we
simplified the implementation to primarily use Level~3 BLAS 
without using memory-bound operations such as pivoting. This made 
our results higher than many available factorization codes thus 
maximizing the influence of the matrix product implemented
with INT8 splitting emulation.
Fig.~\ref{fig:ozkperf1} shows results for the blocking factors $256$, $512$, and $768$.
Fig.~\ref{fig:ozkperf2} shows results for the blocking factors $1024$, $1288$, and $1536$.
Fig.~\ref{fig:ozkperf3} shows results for the blocking factor $2048$. Some of the
blocking factors have prime factors different than $2$ and thus they
divide fewer matrix sizes, which is reflected in the plots with less
measurement points. This is especially visible for blocking factor
$1536$ in Fig.~\ref{fig:ozkperf2}.

\section{Conclusions and Future Work}
\label{sec:conclude}

We presented the numerical and performance results from running HPL's LU
factorization using hardware-native FP64 and INT8 splitting emulation.
We showed how current random number generator poses little problems for
the split integer emulation and we proposed an alternative data
generator that is capable of posing more challenging numerical range of
entries. This resulted in a greater potential of numerical failures of
the emulation and may be indicative of issues possible in applications
using similar dense solvers.

We envision a number of extensions of our work including additional
solver types, matrix properties~\cite{zielke1974maxkond}, and data
formats in need of input data
triggering numerical inaccuracies.

\section*{Acknowledgments}

The authors wish to acknowledge the following individuals for their
contributions and support: Koley Borchard, Chris Berardi, Bob Bond, Alan
Edelman, Peter Fisher, Jeff Gottschalk, Chris Hill, Charles Leiserson,
Kirsten Malvey, Sanjeev Mohindra, Heidi Perry, Christian Prothmann,
Steve Rejto, Scott Ruppel, Daniela Rus, Mark Sherman, Marc Zissman.

\bibliographystyle{unsrt}
\bibliography{mxp}

\begin{thebibliography}{10}

\bibitem{doi:10.1177/10943420211003313}
Ahmad Abdelfattah, Hartwig Anzt, Erik~G Boman, Erin Carson, Terry Cojean, Jack
  Dongarra, Alyson Fox, Mark Gates, Nicholas~J Higham, Xiaoye~S Li, Jennifer
  Loe, Piotr Luszczek, Srikara Pranesh, Siva Rajamanickam, Tobias Ribizel,
  Barry~F Smith, Kasia Swirydowicz, Stephen Thomas, Stanimire Tomov, Yaohung~M
  Tsai, and Ulrike~Meier Yang.
\newblock A survey of numerical linear algebra methods utilizing
  mixed-precision arithmetic.
\newblock {\em The International Journal of High Performance Computing
  Applications}, 2021.

\bibitem{reuther2019llscsuvey}
Albert Reuther, Peter Michaleas, Michael Jones, Vijay Gadepally, Siddharth
  Samsi, and Jeremy Kepner.
\newblock Survey and benchmarking of machine learning accelerators.
\newblock In {\em 2019 IEEE High Performance Extreme Computing Conference
  (HPEC)}, pages 1--9, 2019.

\bibitem{ozaki2012errorfree}
Katsuhisa Ozaki, Takeshi Ogita, Shin'ichi Oishi, and Siegfried~M. Rump.
\newblock Error-free transformations of matrix multiplication by using fast
  routines of matrix multiplication and its applications.
\newblock {\em The Numerical Algorithms}, 59(1):95--118, 2012.

\bibitem{mukunoki2020dgemmtc}
Daichi Mukunoki, Katsuhisa Ozaki, Takeshi Ogita, and Toshiyuki Imamura.
\newblock {DGEMM} using tensor cores, and its accurate and reproducible
  versions.
\newblock In Ponnuswamy Sadayappan, Bradford~L. Chamberlain, Guido Juckeland,
  and Hatem Ltaief, editors, {\em High Performance Computing}, pages 230--248.
  Springer International Publishing, Cham, 2020.

\bibitem{ootomo2022fp32}
Hiroyuki Ootomo and Rio Yokota.
\newblock Recovering single precision accuracy from tensor cores while
  surpassing the {FP32} theoretical peak performance.
\newblock {\em The International Journal of High Performance Computing
  Applications}, 36(4):475--491, 2022.

\bibitem{ozaki2013generrorfree}
Katsuhisa Ozaki, Takeshi Ogita, Shin'ichi Oishi, and Siegfried~M Rump.
\newblock Generalization of error-free transformation for matrix multiplication
  and its application.
\newblock {\em Nonlinear Theory and Its Applications}, 4(1):2--11, 2013.

\bibitem{ozaki2015imperrorfree}
Katsuhisa Ozaki, Takeshi Ogita, and Shin'ichi Oishi.
\newblock Improvement of error-free splitting for accurate matrix
  multiplication.
\newblock {\em Journal of computational and applied mathematics},
  288:127–140, 2015.

\bibitem{ozaki2016errorfreevalid}
Katsuhisa Ozaki, Takeshi Ogita, and Shin'ichi Oishi.
\newblock Error-free transformation of matrix multiplication with a posteriori
  validation.
\newblock {\em Numerical Linear Algebra with Applications}, 23(5):931–946,
  2016.

\bibitem{ichimura2018spsmatvec}
Shuntaro Ichimura, Takahiro Katagiri, Katsuhisa Ozaki, Takeshi Ogita, and Toru
  Nagai.
\newblock Threaded accurate matrix-matrix multiplications with sparse
  matrix-vector multiplications.
\newblock In {\em 2018 IEEE International Parallel and Distributed Processing
  Symposium Workshops (IPDPSW)}, page 1093–1102, 2018.

\bibitem{mukunoki2021matmul128}
Daichi Mukunoki, Katsuhisa Ozaki, Takeshi Ogita, and Toshiyuki Imamura.
\newblock Accurate matrix multiplication on binary128 format accelerated by
  ozaki scheme.
\newblock In {\em Proceedings of the 50th International Conference on Parallel
  Processing (Lemont, IL, USA) (ICPP '21)}, New York, NY, USA, 2021.
  Association for Computing Machinery.
\newblock Article 78, 11 pages.

\bibitem{lange2022altozaki}
Marko Lange and Siegfried~M Rump.
\newblock An alternative approach to ozaki's scheme for error-free
  transformation of matrix multiplication.
\newblock In {\em International Workshop on Reliable Computing and
  Computer-Assisted Proofs}, 2022.

\bibitem{ozaki2022errapps}
Katsuhisa Ozaki, Daichi Mukunoki, and Takeshi Ogita.
\newblock Error handling and applications of iterative error-free
  transformations for matrix multiplication (in {Japanese}).
\newblock In {\em The 18th Joint Meeting of JSIAM Activity Groups}, 2022.

\bibitem{mukunoki2022infprecdot}
Daichi Mukunoki, Katsuhisa Ozaki, Takeshi Ogita, and Toshiyuki Imamura.
\newblock Infinite-precision inner product and sparse matrix-vector
  multiplication using ozaki scheme with {DOT2} on manycore processors.
\newblock In {\em International Conference on Parallel Processing and Applied
  Mathematics}, page 40–54. Springer, 2022.

\bibitem{uchino2024eigsvd}
Yuki Uchino.
\newblock {\em Study on High-Reliability Numerical Methods for
  Eigenvalue/Singular Value Decomposition and Matrix Multiplication (in
  {Japanese})}.
\newblock PhD thesis, Shibaura Institute of Technology, 2024.

\bibitem{ootomo2024dgemmint}
Hiroyuki Ootomo, Katsuhisa Ozaki, and Rio Yokota.
\newblock {DGEMM} on integer matrix multiplication unit.
\newblock {\em The International Journal of High Performance Computing
  Applications}, 38(4):297--313, 2024.

\bibitem{uchino2025ozakiperf}
Yuki Uchino, Katsuhisa Ozaki, and Toshiyuki Imamura.
\newblock Performance enhancement of the {Ozaki} scheme on integer matrix
  multiplication unit.
\newblock {\em The International Journal of High Performance Computing
  Applications}, 39(3):462--476, 2025.

\bibitem{ozaki2025errfreeext}
Katsuhisa Ozaki, Daichi Mukunoki, and Takeshi Ogita.
\newblock Extension of accurate numerical algorithms for matrix multiplication
  based on error-free transformation.
\newblock {\em Japan Journal of Industrial and Applied Mathematics},
  42(1):1--20, 2025.

\bibitem{Dongarra2003}
Jack~J. Dongarra, Piotr Luszczek, and Antoine Petitet.
\newblock The {LINPACK} benchmark: Past, present, and future.
\newblock {\em Concurrency and Computation: Practice and Experience},
  15(9):803--820, August 10 2003.
\newblock {DOI: 10.1002/cpe.728}.

\bibitem{wilkinsonRoundingErrorsAlgebraic1963}
James~H. Wilkinson.
\newblock {\em Rounding Errors in Algebraic Processes}.
\newblock {Prentice-Hall}, {Princeton, NJ, USA}, 1963.

\bibitem{wilkinsonAlgebraicEigenvalueProblem1965}
James~H. Wilkinson.
\newblock {\em The Algebraic Eigenvalue Problem}.
\newblock {Oxford University Press}, {London, UK}, 1965.

\bibitem{trefethenAveragecaseStabilityGaussian1990}
Lloyd~N. Trefethen and Robert~S. Schreiber.
\newblock Average-case stability of {{Gaussian}} elimination.
\newblock {\em SIAM Journal on Matrix Analysis and Applications},
  11(3):335--360, July 1990.

\bibitem{turing1948rounderrmtx}
Alan~M. Turing.
\newblock Rounding-off errors in matrix processes.
\newblock {\em The Quarterly Journal of Mechanics and Applied Mathematics},
  1(1):287--308, January 1948.

\bibitem{fasi2021nopivlugen}
Massimiliano Fasi and Nicholas~J. Higham.
\newblock Matrices with tunable infinity-norm condition number and no need for
  pivoting in lu factorization.
\newblock {\em SIAM Journal on Matrix Analysis and Applications},
  42(1):417--435, 2021.

\bibitem{luszczek2020scaldatagen}
Piotr Luszczek, Yaohung Tsai, Neil Lindquist, Hartwig Anzt, and Jack Dongarra.
\newblock Scalable data generation for evaluating mixed-precision.
\newblock In {\em Proceedings of IEEE HPEC'20: High Performance Extreme
  Computing}, Waltham, MA, USA, 2020.

\bibitem{reuther2018interactive}
A.~{Reuther}, J.~{Kepner}, C.~{Byun}, S.~{Samsi}, W.~{Arcand}, D.~{Bestor},
  B.~{Bergeron}, V.~{Gadepally}, M.~{Houle}, M.~{Hubbell}, M.~{Jones},
  A.~{Klein}, L.~{Milechin}, J.~{Mullen}, A.~{Prout}, A.~{Rosa}, C.~{Yee}, and
  P.~{Michaleas}.
\newblock Interactive supercomputing on 40,000 cores for machine learning and
  data analysis.
\newblock In {\em 2018 IEEE High Performance extreme Computing Conference
  (HPEC)}, pages 1--6, Sep. 2018.

\bibitem{zielke1974maxkond}
G.~Zielke.
\newblock {Testmatrizen mit maximaler Konditionszahl}.
\newblock {\em Computing}, 13(1):33--54, March 1974.

\end{thebibliography}

\end{document}